\documentclass{article}
\usepackage{graphicx}

\title{\Large \bf  Roughness in Cayley Graphs}
\author{ {\bf M.H. Shahzamanian}$^1$, {\bf
M. Shirmohammadi$^2$} and {\bf B. Davvaz}$^1$\\ \\
$^1$Department of Mathematics, Yazd University, Yazd, Iran\\
shahzamanian@stu.yazduni.ac.ir \\
davvaz@yazduni.ac.ir\\
$^2$Department of Computer Engineering, Yazd University, Yazd,
Iran\\
shirmohammadi@stu.yazduni.ac.ir  }

\date{ }
\begin{document}
\maketitle
\begin{abstract}
 In this paper, rough approximations
of Cayley graphs are studied and rough edge Cayley graphs are
introduced. Furthermore, a new algebraic definition called
pseudo-Cayley graphs  containing Cayley graphs is proposed. Rough
approximation is expanded to  pseudo-Cayley graphs. Also, rough
vertex pseudo-Cayley graphs and rough pseudo-Cayley graphs are
introduced. Some theorems are provided, form which some properties
such as connectivity and optimal connectivity are derived. This
approach  opens a new research field in sciences such as data
networks.
\end{abstract}
\section{Introduction}

Graph theory is rapidly moving into the mainstream of mathematics
 mainly because of its applications in diverse fields which include
 biochemistry (genomics), electrical engineering (communications networks and
 coding theory), computer science (algorithms and computations) and
 operations research (scheduling). The wide scope of these and other
 applications has been well-documented cf. \cite{Caccetta, Roberts}. The powerful
 combinatorial methods found in graph theory have also been used to
 prove significant and well-known results in a variety of areas in
 mathematics itself.
In mathematics, the Cayley graph, also known as the Cayley color
graph,
 is the graph that encodes the structure of a discrete group. Its
 definition is suggested by Cayley's theorem (named after Arthur Cayley) and
 uses a particular, usually finite, set of generators for the group. It
 is a central tool in combinatorial and geometric group theory.
\\
\indent The concept of rough set was originally proposed by Pawlak
 \cite{Paw1} as a formal tool
for modelling and processing in complete information in
information
 systems. Since
then the subject has been investigated in many papers (for
example, see
 \cite{D, Deng, Kondo, Mi, Paw2, Pei, Wu, Yang1, Yang2, Yao, Zhang,
 Zhu1, Zhu2}).
The theory of rough set is an extension of set theory, in which a
 subset of a universe is described by a pair of ordinary sets called the
 lower and upper approximations.
A key concept in Pawlak rough set model is an equivalence
relation. The
 equivalence
classes are the building blocks for the construction of the lower
and upper approximations. The lower approximation of a given set
is the
 union
of all the equivalence classes which are subsets of the set, and
the
 upper
approximation is the union of all the equivalence classes which
have a
 nonempty
intersection with the set. It is well known that a partition
induces an equivalence relation on a set and vice versa. The
properties of rough
 sets can
be thus  examined via either partition or equivalence classes. The
 objects
of the given universe $U$ can be divided into three classes with
 respect to any
subset $A\subseteq U$
\begin{itemize}
\item[(1)] the objects, which are definitely in $A$; \item[(2)]
the objects, which are definitely not in $A$; \item[(3)] the
objects, which are possibly in $A$.
\end{itemize}
The objects in class 1 form the lower approximation of $A$, and
the
 objects in
types 1 and 3 together form its upper approximation. The boundary
of
 $A$ contains
 objects in class 3. Rough sets are a suitable mathematical model
of vague concepts, i.e. concepts without sharp boundaries. Rough
set
 theory is
emerging as a powerful theory dealing with imperfect data. It is
an
 expanding
research area which stimulates explorations on both real-world
 applications
and on the theory itself. It has found practical applications in
many areas such as knowledge discovery, machine learning, data
analysis,
 approximate
classification, conflict analysis, and so on. \\
\indent Biswas and Nanda \cite{Biswas} introduced the rough
subgroup
 notion; Kuroki \cite{Ku1} defined the rough ideal in a semigroup; Kuroki
and Wang \cite{rough-group} studied the lower and upper
approximations
 with respect to normal subgroups.
In \cite{rough-ring, D3}, Davvaz concerned a relationship between
rough sets and ring theory and considered a ring as a universal
set and
 introduced the notion of rough ideals and rough subrings
with respect to an ideal of a ring. In \cite{KD}, Kazanc and
Davvaz
 introduced the notions of rough prime (primary) ideals
and rough fuzzy prime (primary) ideals in a ring and gave some
 properties of such ideals. Rough modules have been investigated
by Davvaz and Mahdavipour \cite{DM}. In \cite{Xiao}, the notions
of
 rough prime ideals and rough
fuzzy prime ideals in a semigroup were introduced. Jun \cite{J2}
 discussed the roughness of $\Gamma$-subsemigroups and ideals in
 $\Gamma$-semigroups.
In \cite{J1}, as a generalization of ideals in BCK-algebras, the
notion
 of rough ideals is discussed. In \cite{LD}, Leoreanu-Fotea and Davvaz
 introduced the concept of $n$-ary subpolygroups. For more information
 about algebraic properties of rough sets refer to \cite{D1, D2,
 D4, D5, Deg, KYD, Li, Lin}.\\
\indent   In this paper, rough approximations of Cayley graphs are
studied and rough edge Cayley graphs are introduced. Furthermore,
a new algebraic definition called pseudo-Cayley graphs  containing
Cayley graphs is proposed. Rough approximation is expanded to
pseudo-Cayley graphs. Also, rough vertex pseudo-Cayley graphs and
rough pseudo-Cayley graphs are introduced.
\\
\indent In distributed systems, reliability and fault tolerance
are major factors which have been received considerable attentions
in scientific literatures \cite{hey,joy}. In special cases, data
networks use Cayley graphs in their backbone, concentrating on
edge and vertex connectivity. The vertex connectivity (edge
connectivity),
 is the minimum number of vertices (edges), that must be removed in
order to disconnect the graph. The fault tolerance of a connected
graph is the maximum number k such that, if any k vertices are
removed, the resulting subgraph is still connected. Reliability
focuses on probabilistic edge connectivity. By computing the rough
edge Cayley graphs of a modelled network, some parameters can be
derived concerning edge connectivity. Also, by computing the
vertex rough pseudo-Cayley graphs of a modelled networked, some
parameters
 can be derived with respect to vertex connectivity.


\section{Basic facts about Cayley graphs}

A \textit{graph} is a pair $X = (V(X),E(X))$ of
 sets satisfying $E(X) \subseteq [V(X)]^{2}$; thus, the elements of
 $E(X)$
 are 2-element subsets of $V(X)$. The
 elements of $V(X)$ are \textit{vertices} (or nodes) of the graph
 $X$, and the
 elements of $E(X)$ are its \textit{edges}.
A graph $Y$ is a \textit{subgraph} of $X$ (written $Y\subseteq X$)
if $V(Y) \subseteq V(X)$, $E(Y) \subseteq E(X)$.
 When $Y \subseteq X$ but $Y \neq X$, we write $Y \subset X$ and call
 $Y$
a proper subgraph of $X$. If $Y$ is a subgraph of $X$, $X$ is a
supergraph of $Y$. A spanning subgraph (or spanning supergraph) of
$X$ is a subgraph (or supergraph) $Y$ with $V(Y) = V(X)$.

The \textit{union} $X_{1} \cup X_{2}$ of $X_{1}$ and $X_{2}$ is
the supergraph with vertex set $V(X_{1})
 \cup V(X_{2})$ and edge set $E(X_{1}) \cup E(X_{2})$. The
\textit{intersection} $X_{1} \cap X_{2}$ of $X_{1}$ and $X_{2}$ is
defined similarly, but in this case $X_{1}$ and $X_{2}$ must have
at least one vertex in common.

A \textit{walk} (of length $k$) in a graph $X$ is a non-empty
alternating sequence $v_{0}e_{0}v_{1}e_{1} \ldots e_{k-1}v_{k}$ of
vertices and edges in $X$ such that $e_{i} = (v_{i},v_{i+1})$ for
all $i<k$. If $v_{0} = v_{k}$, the walk is closed. If the vertices
in a walk are all distinct, it defines an obvious \textit{path} in
$X$. In general, every walk between two vertices contains a path
between these vertices.

A non-empty graph $X$ is called \textit{connected} if any two of
its vertices are linked by a path in $X$. A connected graph $X$ is
called \textit{optimal connected} if every spanning subgraph of
$X$ is not connected.\\ \\
{\bf Definition 2.1.} Taking any finite group $G$, let $S \subset
G$ be such that $1 \not \in S$ (where $1$ represents the identity
element of $G$) and $s \in S$ implies that $s^{-1} \in S$ (where
$s^{-1}$ represents the inverse element of $s$). The
\textit{Cayley graph} $(G;S)$ is a graph whose vertices are
labelled with the elements of $G$, in which there is an edge
between two vertices $g$ and $gs$  if and only if $s \in S$.\\
\indent The exclusion of $1$ from $S$ eliminates the possibility
of loops in the graph. The inclusion of the inverse of any element
which is itself in $S$ means that an edge is in the graph
regardless of which end vertex is considered.\\
\indent Let $R$ be not group and be a subset of $G$, if $R$
contains $S$ and $SR \subseteq R$ where $SR=\{sr |s \in S, r \in
R\}$, then the \textit{pseudo-Cayley graph} $(R;S)$ is a graph
whose vertices are labelled with the elements of $R$, in which
there is an edge between two vertices $r$ and $rs$ if and only if
$s \in S$.
\\ \\
\indent Let $G$ be a group and $X$ be a subset of $G$. Let
$\{H_{i} \mid i \in I\}$ be the family of all subgroups of $G$
which contains $X$. Then $\bigcap_{i \in I} H_{i}$ is called the
\emph{subgroup} of $G$ \emph{generated by the set} $X$ and denoted
$<X>$. Obviously, if $<X>=G$, then $X$ generates $G$.
\\ \\
{\bf Theorem 2.2.}\cite{godsil}. A Cayley graph $(G;S)$ is
connected if and only if $S$ generates $G$.\\ \\
\indent  A subset $S$ of $G$ is called \textit{minimal Cayley set}
if it generates $G$, and $S\backslash \{ s , s^{-1} \}$ generates
a proper subgroup of $G$ for all $s \in S$.\\ \\
{\bf Theorem 2.3.}\cite{joy}. If $S$ is minimal Cayley set for the
finite group $G$, then the Cayley graph $(G;S)$ has optimal
connectivity.\\ \\
{\bf Theorem 2.4.} If $X_{1}=(G;S_{1})$ and $X_{2}=(G;S_{2})$ are
Cayley graphs, then
\begin{itemize}
    \item[(1)] $X_{1} \cup X_{2} = (G;S_{1} \cup S_{2})$,
    \item[(2)] $X_{1} \cap X_{2} = (G;S_{1} \cap S_{2})$.
\end{itemize}
{\it  Proof.} (1) Let $e$ be an edge of $(G;S_{1} \cup S_{2})$
then there exist $g \in G$, and $s \in S_{1} \cup S_{2}$ such that
$e$ is connecting two vertices $g$ and $gs$. Since $s \in  S_{1}
\cup S_{2}$ then $s \in S_{1}$ or $s \in S_{2}$, which  equals by
$e \in E(X_{1})$ or $e \in E(X_{2})$. Therefore according to
definition, $e \in E(X_{1} \cup X_{2})$. Conversely, in the
similar way, any edge of $E(X_{1} \cup X_{2})$ is an edge of
$(G;S_{1} \cup S_{2})$. This result and $V(X_{1} \cup X_{2}) =
V((G;S_{1} \cup S_{2})) = G$, yield $X_{1} \cup X_{2} = (G;S_{1}
\cup S_{2})$.
\\
\indent (2) It is straightforward.\\
 \indent Notice that $X_{1}
\cup X_{2}$ and $X_{1} \cap X_{2}$ are Cayley
graphs.\\ \\
{\bf Theorem 2.5.} If $X_{1}=(H_{1};S)$ and $X_{2}=(H_{2};S)$
$(H_{1} , H_{2} \leq G$ which means $H_{1}$ and $H_{2}$ are
subgroups of $G$) are Cayley graphs, then
\begin{itemize}
    \item[(1)] $X_{1} \cup X_{2} = (H_{1} \cup H_{2};S)$,
    \item[(2)] $X_{1} \cap X_{2} = (H_{1} \cap H_{2};S)$.
\end{itemize}
{\it  Proof.} The proof is similar to 2.4.

Notice that $X_{1} \cap X_{2}$ is a Cayley graph, but $X_{1} \cup
X_{2}$ may not be a Cayley graph. $X_{1} \cup X_{2}$ always is a
pseudo-Cayley graph.\\ \\
{\bf Theorem 2.6.} If $X_{1}=(H_{1};S_{1})$ and
$X_{2}=(H_{2};S_{2})$ $(H_{1} , H_{2} \leq G)$ are Cayley graphs,
then $X_{1} \cap X_{2} = (H_{1} \cap H_{2};S_{1} \cap S_{2})$.\\
{\it  Proof.} The proof is similar to 2.4.\\
\\
\indent Notice that $X_{1} \cap X_{2}$ is a Cayley graph. $X_{1}
\cup
X_{2}$ may not be a pseudo-Cayley graph.\\ \\
{\bf Theorem 2.7.} If $X_{1}=(G;S_{1})$, $X_{2}=(G;S_{2})$,
$Y_{1}=(G_{1};S)$ and $Y_{1}=(G_{2};S)$ are Cayley graphs, then
\begin{itemize}
 \item[(1)] $X_{1} \subseteq X_{2}$  if and only if $S_{1}\subseteq
 S_{2}$,
 \item[(2)] $Y_{1} \subseteq Y_{2}$  if and only if $G_{1} \subseteq
 G_{2}$.
\end{itemize}
{\it  Proof.} (1) Let $S_{1}$ be a subset of $S_{2}~(S_{1}
\subseteq S_{2})$. Suppose that $e$ is an arbitrary  edge of
$E(X_{1})$, then there exist $g \in G$ and $s_{1} \in S_{1}$ such
that $e = (g;gs_{1})$. Since $s_{1} \in S_{1} \subseteq S_{2}$,
then $e \in E(X_{2})$. Therefore $E(X_{1})$ will be a subset of
$E(X_{2})$.

Conversely, let $E(X_{1})$ be a subset of $E(X_{2})~(E(X_{1})
\subseteq E(X_{2}) )$. Suppose that $s_{1}$ is any element of
$S_{1}$. For every $g \in G$, we have $(g;gs_{1}) \in E(X_{1})$.
Therefore this gives $(g;gs_{1}) \in E(X_{2})$ and as a result
$s_{1} \in S_{2}$ and then $S_{1} \subseteq S_{2}$. This result
and $V(X_{1}) = V(X_{2}) = G $ lead $X_{1} \subseteq X_{2}$ if and
only if
$S_{1} \subseteq S_{2}$.\\
\indent (2) It is straightforward.\\
\section{Rough groups}

If $N$ is a subgroup of a group $G$, then the following conditions
are equivalent.
\begin{itemize}
\item[(1)]  Left and right congruence modulo $N$ coincide (that
is,
    define the same equivalence relation on $G$),
\item[(2)]  $aN = Na$ for all $a \in G$, \item[(3)] for all $a \in
G$, $aNa^{-1} \subseteq N$, where $aNa^{-1} = \{ana^{-1}|n \in
    N\}$,
\item[(4)] for all $a \in G$, $aNa^{-1} = N$.
\end{itemize}
A subgroup $N$ of a group $G$ which satisfies the above equivalent
conditions  is said to be \textit{normal} in $G$ (or a
\textit{normal subgroup} of $G$). Let $H$ and $N$ be normal
subgroups of a group $G$. Then as it is well known and easily
seen, $H \cap N$ is also a normal subgroup of $G$.
\\ \\
\indent  Let $G$ be a group (as universe) with identity $1$ and
$N$ be a normal subgroup of $G$. If $A$ is a nonempty subset of
$G$, then the sets $$N_{-}(A) = \{x \in G \ | \  xN \subseteq A \}
\ \ {\rm  and} \ \  N^{\wedge} (A) = \{ x \in G \ | \  xN \cap A
\neq \emptyset \}$$ are called, respectively, {\it lower} and {\it
upper approximations of a set $A$ with respect to the normal
subgroup $N$}.\\ \\
{\bf Theorem 3.1.} \cite{rough-group}. Let $H$ and $N$ be normal
subgroups of a group $G$. Let $A$ and $B$ be any nonempty subsets
of $G$. Then
\begin{itemize}
    \item[(1)] $N_{-}(A) \subseteq A \subseteq N^{\wedge} (A)$,
    \item[(2)] $N^{\wedge} (A \cup B) = N^{\wedge} (A) \cup N^{\wedge}
    (B)$,
    \item[(3)] $N_{-}(A \cap B) = N_{-}(A) \cap N_{-}(B)$,
    \item[(4)] $A \subseteq B$ implies $N_{-}(A) \subseteq N_{-}(B)$,
    \item[(5)] $A \subseteq B$ implies $N^{\wedge} (A) \subseteq
 N^{\wedge}
    (B)$,
    \item[(6)] $N_{-}(A \cup B) \supseteq N_{-}(A) \cup
    N_{-}(B)$,
    \item[(7)]  $N^{\wedge} (A \cap B) \subseteq N^{\wedge} (A)
    \bigcap N^{\wedge} (B)$,
    \item[(8)]  $N \subseteq H $ implies $ N^{\wedge} (A) \subseteq
    H^{\wedge} (A)$,
    \item[(9)] $N \subseteq H $ implies $H_{-}(A) \subseteq N_{-}(A)$.
\end{itemize}

\noindent {\bf Theorem 3.2.}\cite{rough-group}. Let $H$ and $N$ be
normal subgroups of a group $G$. If $A$ is a non-empty subset of
$G$, then
\begin{itemize}
    \item[(1)] $(H \cap N)^{\wedge} (A) = H^{\wedge} (A) \cap
 N^{\wedge}
    (A)$,
    \item[(2)] $(H \cap N)_{-}(A) = H_{-}(A) \cap N_{-}(A)$.
\end{itemize}

$N(A) = (N_{-}(A),N^{\wedge} (A))$ is called a {\it rough set} of
$A$ in $G$. A non-empty subset $A$ of a group $G$ is called an
{\it $N^{\wedge}$-rough (normal) subgroup} of $G$ if the upper
approximation  of $A$ is a (normal) subgroup of $G$. Similarly, a
nonempty subset $A$ of $G$ is called an {\it $N_{-}$-rough
(normal) subgroup} of $G$ if lower approximation is a (normal)
subgroup of $G$.\\ \\
{\bf Theorem 3.3.}\cite{rough-group}. Let  $N$ be a normal
subgroup of a group $G$.
\begin{itemize}
\item[(1)] If $A$ is a subgroup of $G$, then it is an
$N^{\wedge}$rough subgroup of $G$. \item[(2)] If $A$ is a normal
subgroup of $G$, then it is an $N^{\wedge}$-rough normal subgroup
of $G$. \item[(3)] If $A$ is a subgroup of $G$ such that
$N\subseteq A$, then it is an $N_{-}$-rough subgroup of $G$.
\item[(4)] If $A$ is a normal subgroup of $G$ such that
$N\subseteq A$, then it is an $N_{-}$-rough normal subgroup of
$G$.
\end{itemize}
\section{Rough edge Cayley graphs}
\begin{figure}
\begin{center}
\includegraphics[1mm,0mm][150mm,30mm]{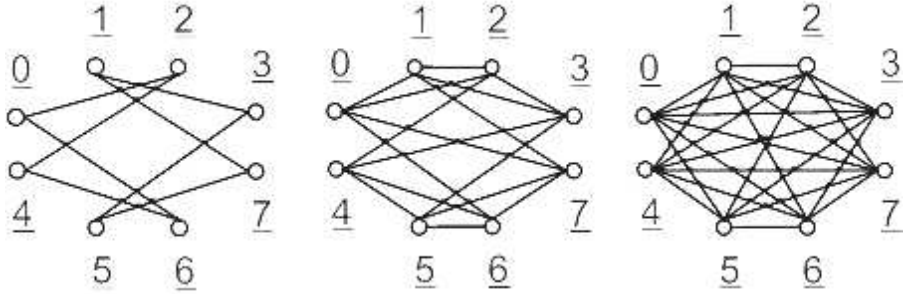}
\caption{The above  graphs are, respectively $\underline{X}$, $X$
and $\overline{X}$.} \label{fig1}
\end{center}
\end{figure}

In this section, concept of lower and upper approximations edge
Cayley graphs of a Cayley graph with respect to a normal subgroup
is discussed  then some
 properties of the lower and upper approximations are brought.\\ \\
{\bf Definition 4.1.} Let $G$ be a finite group with identity $1$,
$N$ be a normal subgroup of $G$ and $X=(G;S)$ be a Cayley graph.
Then the following graphs (we will prove these graphs are Cayley
graphs)
$$\overline{X} = (G;N^{\wedge} (S)^{\ast})\ \  {\rm  and} \ \ \underline{X} =
(G;N_{-}(S))$$ where $N^{\wedge} (S)^{\ast} = N^{\wedge}
(S)\backslash 1$, are called, respectively,  {\it lower and upper
approximations edge Cayley graphs} of the Cayley graph $X(G;S)$
with respect to the normal subgroup $N$.\\ \\
{\bf Theorem 4.2.} The two graphs $\underline{X}$  and
$\overline{X} $
are Cayley graphs.\\ \\
{\it Proof.} Theorem 3.1 (1) gives $N_{-}(S) \subseteq S$. Then $1
\not \in N_{-}(S)$. Suppose that $s$ is an arbitrary element of
$N_{-}(S)$. Then $sN \subseteq S $ which implies that $sn^{-1} \in
S$, for all $n \in N$. Then $(sn^{-1})^{-1} = ns^{-1} \in S$ and
so $Ns^{-1} \subseteq S$ or $s^{-1}N \subseteq S$. Therefore
$s^{-1} \in N_{-}(S)$.
\\
\indent Now, suppose that $s$ is an arbitrary element of
$N^{\wedge} (S)^{\ast}$. Then $sN \cap S \neq \emptyset$ which
implies that there exists $a \in sN \cap S$. Hence there exists $n
\in N$ such that $a = sn \in S$, so $a^{-1} = n^{-1}s^{-1} \in S$.
On the other hand, we have $n^{-1}s^{-1} \in Ns^{-1} = s^{-1}N$.
Thus $n^{-1}s^{-1} \in S \cap s^{-1}N $, which implies that
$s^{-1}N \cap S \neq \emptyset$ and so $s^{-1} \in N^{\wedge}
(S)^{\ast}$. \\
\indent Therefore,  $\underline{X} $ and $\overline{X}$ are Cayley
 graphs.\\ \\
{\bf Example 4.3.} Let $G$ be a group congruence modulo $8$
integral number ${Z}$. Let $N = \{\underline{0} , \underline{4} \}
$ be a normal subgroup of $G$ and Cayley graph $X = (G;S)$ such
that $S$ equals to $\{ \underline{1} , \underline{2} ,
\underline{6} , \underline{7} \}$. We have $\overline{X} =
(G;\{\underline{1} , \underline{2} , \underline{3} , \underline{5}
, \underline{6} , \underline{7} \})$ and $ \underline{X} = (G; \{
\underline{2} ,\underline{6}\})$ (See figure \ref{fig1}).

\noindent {\bf Theorem 4.4.} Let $N$ and $H$ be normal subgroups
of a group $G$. Let $X=(G;S)$, $X_{1}=(G;S_{1})$ and
$X_{2}=(G;S_{2})$ be Cayley graphs. Then we have
\begin{itemize}
    \item[(1)] $\underline{X} \subseteq X \subseteq \overline{X}$,
\item[(2)] $\overline{X_{1} \cup X_{2}} = \overline{X_{1}} \cup
\overline{X_{2}} $,
 \item[(3)] $\underline{X_{1} \cap X_{2}} = \underline{X_{1}} \cap
\underline{X_{2}} $,

 \item[(4)] $X_{1} \subseteq X_{2} \Rightarrow \underline{X_{1}}
 \subseteq \underline{X_{2}}$,

 \item[(5)] $X_{1} \subseteq X_{2} \Rightarrow \overline{X_{1}}
 \subseteq \overline{X_{2}}$,

\item[(6)] $\underline{X_{1} \cup X_{2}} \supseteq
\underline{X_{1}} \cup \underline{X_{2}} $,

\item[(7)] $\overline{X_{1} \cap X_{2}} \subseteq \overline{X_{1}}
\cap \overline{X_{2}} $,

\item[(8)] $N \subseteq H \Rightarrow \overline{X_{N} } \subseteq
\overline{X_{H}}$,

\item[(9)]  $N \subseteq H \Rightarrow \underline{X_{H} }
\subseteq \underline{X_{N}}$.

\end{itemize}
{\it  Proof.}
\begin{itemize}
    \item[(1)] By Theorem 3.1 (1), $N_{-}(S) \subseteq S \subseteq
N^{\wedge} (S)$ then $N_{-}(S) \subseteq S \subseteq N^{\wedge}
(S)^{\ast}$ $(1 \not \in S)$ and Theorem 2.7 (1) leads
$\underline{X} \subseteq X \subseteq \overline{X}$.

\item[(2)] Theorem 2.4 gives $X_{1} \cup X_{2} = (G;N^{\wedge}
(S_{1})^{\ast} \cup N^{\wedge} (S_{2})^{\ast})$. By Theorem 3.1
(5), we have $N^{\wedge} (S_{1})^{\ast}$ and $N^{\wedge}
(S_{2})^{\ast} $ $\subseteq$ $ N^{\wedge} (S_{1} \cup
S_{2})^{\ast}$. Then Theorem 2.7 shows $\overline{X_{1}} \cup
\overline{X_{2}} \subseteq \overline{X_{1} \cup X_{2}}$.
Conversely, according to Theorem 3.1 (2), we have $N^{\wedge}
(S_{1})^{\ast} \cup N^{\wedge} (S_{2})^{\ast} = N^{\wedge} (S_{1}
\cup S_{2})^{\ast}$. Suppose that $(g;gs)$ is any edge of
$E(\overline{X_{1} \cup X_{2}})$ and $ s \in N^{\wedge} (S_{1}
\cup S_{2})^{\ast}$. Then we obtain $s \in N^{\wedge}
(S_{1})^{\ast} \cup N^{\wedge} (S_{2})^{\ast}$ and $s \in
N^{\wedge} (S_{1})^{\ast} $ or $s \in N^{\wedge} (S_{2})^{\ast} $.
Therefore $(g;gs)$ is any edge of $\overline{X_{1}}$ or
$\overline{X_{2}}$. Finally, we have $\overline{X_{1} \cup X_{2}}
= \overline{X_{1}} \cup \overline{X_{2}}$.

\item[(3)]  By Theorem 3.1 (3), the proof is similar to (2).

\item[(4)] Assume that $X_{1} \subseteq X_{2}$. Then $S_{1 }
\subseteq S_{2}$, which implies that $N_{-}(S_{1 }) \subseteq
N_{-}(S_{2})$. Hence $\underline{X_{1}} \subseteq
\underline{X_{2}}$.

\item[(5)] By Theorem 3.1 (5), the proof is similar to (4).

\item[(6)]  Theorem 3.1 (6) gives  $N_{-}(S_{1 }) \cup
N_{-}(S_{2}) \subseteq N_{-}(S_{1 } \cup S_{2})$. Then $N_{-}(S_{1
}) \subseteq N_{-}(S_{1 } \cup S_{2})$ and $N_{-}(S_{2}) \subseteq
N_{-}(S_{1 } \cup S_{2})$, therefore we have $\underline{X_{1}
\cup X_{2}} \supseteq \underline{X_{1}}$ and $\underline{X_{1}
\cup X_{2}} \supseteq \underline{X_{2}}$. And finally,
$\underline{X_{1} \cup X_{2}} \supseteq \underline{X_{1}} \cup
\underline{X_{2}}$.

\item[(7)]  By Theorem 3.1 (7),  the proof is similar to (6).

\item[(8)] Assume that $N \subseteq H$.  Theorem 3.1 (8) yields
$N^{\wedge}(S) \subseteq H^{\wedge}(S)$. Then
$N^{\wedge}(S)^{\ast} \subseteq H^{\wedge}(S)^{\ast}$. Now, based
on Theorem 2.7 (1), we obtain $\overline{X_{N} } \subseteq
\overline{X_{H}}$.

\item[(9)] By Theorem 3.1 (9), the proof is similar to (8).
\end{itemize}

\noindent {\bf Example 4.5.} Here, we present some examples which
show
 the contradiction of
the converse part of the above items(4--9). Let $G$ be a dihedral
group with order 6.
\begin{itemize}

 \item[(1)]  $X_{1} \not \subseteq X_{2} , \underline{X_{1}}
  \subseteq \underline{X_{2}}:~
  X_{1}=(G;\{\varepsilon\}),~X_{2}=(G;\{P\varepsilon\}),~N=\{1,P,P^{2}\}.$
 \item[(2)]  $X_{1} \not \subseteq X_{2} , \overline{X_{1}}
 \subseteq \overline{X_{2}}:~
  X_{1}=(G;\{\varepsilon\}),~X_{2}=(G;\{P\varepsilon\}),~N=G.$
\item[(3)]  $\underline{X_{1}} \cup \underline{X_{2}} \not
\supseteq \underline{X_{1} \cup X_{2}}:~
  X_{1}=(G;\{\varepsilon,P^{2}\varepsilon\}),~X_{2}=(G;\{\varepsilon,P\varepsilon\}),~N=\{1,P,P^{2}\}$.
\item[(4)]  $\overline{X_{1}} \cap \overline{X_{2}} \not \subseteq
\overline{X_{1} \cap X_{2}}:~
  X_{1}=(G;\{P^{2}\varepsilon\}),~X_{2}=(G;\{P\varepsilon\}),~N=\{1,P,P^{2}\}$.
\item[(5)] $N \not \subseteq H , \overline{X_{N} } \subseteq
\overline{X_{H}}:~
  X=(G;\{P,P^{2},\varepsilon,P\varepsilon,P^{2}\varepsilon\}),~N=\{1,P,P^{2}\},~H=\{1\}$.
\item[(6)]  $N \not \subseteq H , \underline{X_{H} } \subseteq
\underline{X_{N}}:~
  X=(G;\{\varepsilon,P\varepsilon,P^{2}\varepsilon\}),~N=\{1,P,P^{2}\},~H=\{1\}$.
\end{itemize}

\noindent {\bf Theorem 4.6.} Let $N$ and $H$ be normal subgroups
of a group $G$. Let $X=(G;S)$ be a Cayley graph. Then
\begin{itemize}
    \item[(1)]  $\overline{X}_{H \cap N} = \overline{X}_{H } \cap
    \overline{X}_{N}$,
    \item[(2)] $\underline{X}_{H \cap N} = \underline{X}_{H} \cap
    \underline{X}_{N}$.
\end{itemize}
{\it Proof.}
\begin{itemize}
    \item[(1)] We have
    $$
    \begin{array}{ll}
    \overline{X}_{H \cap N} & = (G;(H \cap
    N)^{\wedge}(S))\\
    & = (G; H^{\wedge}(G) \cap N^{\wedge}(S))\\
    & = (G;H^{\wedge}(S)) \cap (G;N^{\wedge}(S))\\
    & =\overline{X}_{H} \cap \overline{X}_{N}.
    \end{array}
    $$
    \item[(2)] According to Theorem 3.2 (2), the proof is similar to
 (1).
\end{itemize}

$(\underline{X},\overline{X})$ is called a {\it rough edge Cayley
graph} of $X
 =
(G;S)$. A Cayley graph $X = (G;S)$  is called an {\it
$N^{\wedge}$-edge rough generating} if the $N^{\wedge} (S)^{\ast}$
is a generating set for $G$. Similarly, a Cayley graph $X = (G;S)$
is called an {\it $N_{-}$-edge rough generating} if $N_{-}(S)$ is
a generating set
for $G$.\\
\indent A Cayley graph $X = (G;S)$  is called an {\it
$N^{\wedge}$-edge
 rough
optimal connected} if the $N^{\wedge} (S)^{\ast}$ is a minimal
Cayley set  for $G$. Similarly, a Cayley graph $X = (G;S)$ is
called an {\it $N_{-}$-edge rough optimal} if $N_{-}(S)$ is a
minimal
Cayley set for $G$.\\ \\
{\bf Theorem  4.7.} Let $X=(G;S)$ be a Cayley graph. If $X$ is an
$N^{\wedge}$-edge rough generating, then $\overline{X} $ is
connected. Similarly, If $X$ is an $N_{-}(S)$-edge rough
generating, then $\underline{X} $
is connected.\\ \\
{\it Proof.} It is straightforward.\\ \\
{\bf Theorem 4.8.} Let $X=(G;S)$ be a Cayley graph. If $X$ is an
$N^{\wedge}$-edge rough optimal connected, then $\overline{X} $ is
optimal connected. Similarly, if $X$ is an $N_{-}(S)$-edge rough
optimal connected,
then $\underline{X} $ is optimal connected.\\ \\
{\it Proof.} It is straightforward.
\section{Rough vertex pseudo-Cayley graphs}
In this section, concept of lower and upper approximations vertex
pseudo-Cayley graphs of a pseudo-Cayley graph with respect to a
normal subgroup is introduced. We prove the lower and upper
 approximations are pseudo-Cayley graphs, too. Then  some
 properties of lower and upper approximations are brought.\\ \\
{\bf Definition 5.1.} Let $G$ be a finite group with identity $1$,
$N$ be a normal subgroup, $R$ be a subset of $G$ and $X=(R;S)$ be
a pseudo-Cayley graph. Then the following graphs (we will prove
these
 graphs are pseudo-Cayley
graphs):

$$\overline{X}' = (N^{\wedge} (R);S)\ \ {\rm  and} \ \ \underline{X}' =
(N_{-}(R);S \cap N_{-}(R))$$ are called, respectively,  {\it
lower} and {\it upper approximations vertex pseudo-Cayley graphs}
of a pseudo-Cayley
graph $X(R;S)$ with respect to the normal subgroup $N$.\\ \\
{\bf Theorem 5.2.} The two graphs $\underline{X}'$  and
$\overline{X}' $ are pseudo-Cayley
graphs.\\ \\
{\it Proof.} The definition of pseudo-Cayley graph says  $S
\subseteq R$. Then $S \subseteq N^{\wedge}(R)$. If $ s \in S$ and
$a \in N^{\wedge}(R)$, then $aN \cap R \neq \emptyset$. Hence
there exists $n \in N$ such that $an \in R$. Since $X$ is a
pseudo-Cayley graph then $SR \subseteq R$. So $san \in R $ which
implies that $saN \cap R \neq \emptyset$, so $sa \in
N^{\wedge}(R)$. Thus $SN^{\wedge}(R) \subseteq N^{\wedge}(R)$.
Then $\overline{X}'$ is a pseudo-Cayley graph.

\begin{figure}[h]
\begin{center}
\includegraphics[1mm,0mm][85mm,30mm]{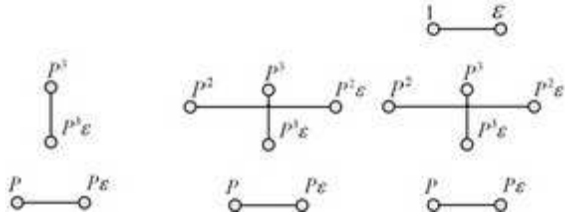}
\caption{The above graphs are, respectively $\underline{X}'$, $X$
and $\overline{X}'$.} \label{fig2}
\end{center}
\end{figure}

If $s \in S$ and $ a \in N_{-}(R)$, then $aN \subseteq R$. So for
all $n \in N $ we have $an \in R$. Since $X$ is a pseudo-Cayley
graph, then $SR \subseteq R$. Thus $san \in R$, which implies that
$saN \subseteq R$ and so $sa \in N_{-}(R)$. Hence $SN_{-}(R)
\subseteq N_{-}(R)$. Therefore $S \cap N_{-}(R)$ has the necessary
conditions for being a
pseudo-Cayley graph. So $\overline{X}'$ is a pseudo-Cayley graph.\\ \\
{\bf Example 5.3.} Let $G$ be a dihedral group with order $8$,  $N
= \{1, P^{2} \} $ be a normal subgroup of $G$ and $R = \{P, P^{2},
P^{3} ,P \varepsilon , P^{2} \varepsilon , P^{3} \varepsilon  \}$
be a subset of $G$. Let pseudo-Cayley graph $X = (R;S)$ such that
$S$ equals to $\{\varepsilon \}$. We have
$$\underline{X}' = (\{P, P^{3} ,P \varepsilon , P^{3} \varepsilon
\};S)$$ and $$ \overline{X}' = (\{P, P^{2}, P^{3}, P \varepsilon ,
P^{2} \varepsilon , P^{3} \varepsilon , 1, \varepsilon \};S)$$
(See figure \ref{fig2}).
\\

In order to find some illustrative examples and contradictions, a
software is developed in C++ programming language. Five classes
include CayleyGraph, PseudoCayleyGraph, Group, NormalSubGroup, and
SubsetS, with their properties and methods are defined. Since
pseudo-Cayley graphs contain Cayley graphs, the CayleyGraph class
inherits PseudoCayleyGraph class. In main thread of program
running mode, the user inputs the number of elements of the group.
The binary operation of the group is then initialized. The groups,
here are restricted to dihedral and congruence groups. Then all
normal subgroups and subsets which satisfy the conditions of
$(G;S)$ to be  Cayley graph, are computed. With the determined
group, and with all computed $N$s and $S$s, the software plots all
possible Cayley graphs. In addition, lower and upper
approximations of the determined group, Rough edge Cayley graph
and  Rough vertex pseudo-Cayley graph can also be computed.
\\

\noindent {\bf Theorem 5.4.} Let $N$ and $H$ be normal subgroups
and $R$, $R_{1}$ and $R_{2}$ be subsets of a group $G$ and $S
\subseteq R, R_{1}, R_{2}$. Let $X=(R;S)$, $X_{1}=(R_{1};S)$ and
$X_{2}=(R_{2};S)$ be pseudo-Cayley graphs. Then we have
\begin{itemize}
    \item[(1)] $\underline{X}' \subseteq X \subseteq \overline{X}'$,
\item[(2)]  $\overline{X_{1} \cup X_{2}}' = \overline{X_{1}}' \cup
\overline{X_{2}}' $,
 \item[(3)]  $\underline{X_{1} \cap X_{2}}' = \underline{X_{1}}' \cap
\underline{X_{2}}' $,

 \item[(4)]  $X_{1} \subseteq X_{2} \Rightarrow \underline{X_{1}}'
 \subseteq \underline{X_{2}}'$,

 \item[(5)]  $X_{1} \subseteq X_{2} \Rightarrow \overline{X_{1}}'
 \subseteq \overline{X_{2}}'$,

\item[(6)]  $\underline{X_{1} \bigcup X_{2}}' \supseteq
\underline{X_{1}}' \bigcup \underline{X_{2}}' $,

\item[(7)]  $\overline{X_{1} \bigcap X_{2}}' \subseteq
 \overline{X_{1}}'
\bigcap \overline{X_{2}}' $,

\item[(8)]  $N \subseteq H \Rightarrow \overline{X_{N} }'
\subseteq \overline{X_{H}}'$,

\item[(9)]  $N \subseteq H \Rightarrow \underline{X_{H} }'
\subseteq \underline{X_{N}}'$.

\end{itemize}
{\it Proof.}
\begin{itemize}
\item [(1)] We have $N_{-}(R) \subseteq R$. Then by Theorem 2.7
(2), we
    have  $(N_{-}(R);N_{-}(R) \cap S) \subseteq (R;N_{-}(R) \cap
    S)$, and Theorem 2.7 (1) yields $(R;N_{-}(R) \cap
    S) \subseteq (R;S)$. So $\underline{X}' \subseteq X$
and Theorem 3.1 (1) shows $R \subseteq N^{\wedge}(R)$. Hence
$(R;S)
 \subseteq (N^{\wedge}(R);S)$.
Therefore $X \subseteq \overline{X}'$.

\item[(2)]   We have
$$
\begin{array}{ll}
\overline{X_{1} \cup X_{2}}' & = (N^{\wedge}(R_{1} \cup
R_{2});S)\\
&= (N^{\wedge}(R_{1}) \cup N^{\wedge}(R_{2});S) \\
& = (N^{\wedge}(R_{1});S) \cup (N^{\wedge}(R_{2});S) \\
&=\overline{X_{1}}' \bigcup \overline{X_{2}}' .
\end{array}
$$

\item[(3)]  By Theorem 3.1 (3), the proof is similar to (2).

\item[(4)] Since  $X_{1} \subseteq X_{2}$, we obtain $R_{1 }
\subseteq
 R_{2}$ (According to Theorem 2.7 (2)) and so $N_{-}(R_{1 }) \subseteq
 N_{-}(R_{2 })$. Then $\underline{X_{1}}' \subseteq \underline{X_{2}}'$.

\item[(5)]  By Theorem 3.1 (5), the proof is similar to (4).

\item[(6)] We have
$$
\begin{array}{ll}
\underline{X_{1} \cup X_{2}}' & = (N_{-}(R_{1 } \cup
R_{2});N_{-}(R_{1 } \cup R_{2}) \cap S)\\
& \supseteq (N_{-}(R_{1}) \cup N_{-}(R_{2 }); (N_{-}(R_{1})
\cup N_{-}(R_{2 })) \cap S )\\
& = (N_{-}(R_{1});(N_{-}(R_{1}) \cup N_{-}(R_{2 })) \cap S) \cup
(N_{-}(R_{2});(N_{-}(R_{1}) \cup N_{-}(R_{2
})) \cap S) \\
& = (N_{-}(R_{1});N_{-}(R_{1}) \cap S) \cup
(N_{-}(R_{2});N_{-}(R_{2}) \bigcap S) \\
& =\underline{X_{1}}' \bigcup \underline{X_{2}}' .
\end{array}
$$

\item[(7)] We have
$$
\begin{array}{ll}
\overline{X_{1} \cap X_{2}}' & = (N^{\wedge}(R_{1})
\cap N^{\wedge}(R_{2}); S)\\
& \subseteq (N^{\wedge}(R_{1});S) \cap (N^{\wedge}(R_{2});S)
\\
& = \overline{X_{1}}' \bigcap \overline{X_{2}}' .
\end{array}
$$

\item[(8)] Since  $N \subseteq H$, then $N^{\wedge}(R) \subseteq
 H^{\wedge}(R)$. Thus
$$\overline{X_{N} }' = (N^{\wedge}(R);S) \subseteq
(H^{\wedge}(R);S) = \overline{X_{H}}'.$$

\item[(9)]  By Theorem 3.1 (9), the proof is similar to (8).

\end{itemize}
\noindent In the following, we present some examples which show
the
 contradiction of
converse parts of the above items(4--9). \\ \\
{\bf Example 5.5.} Let $G=\{1,\varepsilon , P, P \varepsilon ,
P^2, P^2
 \varepsilon , P^3, P^3  \varepsilon \}$ be a dihedral
group with order 8.
\begin{itemize}
\item[(1)] $X_{1} \not \subseteq X_{2} , \underline{X_{1}}'
 \subseteq \underline{X_{2}}':~
  X_{1}=(\{1,P^{2}\};\emptyset),~X_{2}=(\{1,P^{2},\varepsilon\};\emptyset),~N=\{1,P^{2}\}.$
 \item[(2)]  $X_{1} \not \subseteq X_{2} , \overline{X_{1}}'
 \subseteq \overline{X_{2}}':~
  X_{1}=(\{1,\varepsilon,P^{2},P^{2}\varepsilon\};\emptyset),~X_{2}=(\{1,P^{2},P^{2}\varepsilon\};\emptyset),~N=\{1,P^{2}\}.$
\item[(3)]  $\underline{X_{1}}' \bigcup \underline{X_{2}}' \not
\supseteq \underline{X_{1} \bigcup X_{2}}':~
  X_{1}=(\{1,P,P^{2}\};\emptyset),~X_{2}=(\{1,P^{3}\};\emptyset),~N=\{1,P^{2}\}$.
\item[(4)]  $ \overline{X_{1}}' \bigcap \overline{X_{2}}' \not
\subseteq \overline{X_{1} \bigcap X_{2}}':~
  X_{1}=(\{1,P,P^{2}\};\emptyset),~X_{2}=(\{1,P^{2},P^{3}\};\emptyset),~N=\{1,P^{2}\}$.
\item[(5)]  $N \not \subseteq H , \overline{X_{N} }' \subseteq
\overline{X_{H}}':~X=(\{1,P,P^{2},P^{3}\};\emptyset),~H=\{1,P^{2}\},~N=\{1,P,P^{2},P^{3}\}$.
 \item[(6)]  $N \not \subseteq H ,
\underline{X_{H} }' \subseteq \underline{X_{N}}'
:~X=(\{1,P,P^{2},P^{3}\};\emptyset),~H=\{1,P^{2}\},~N=\{1,P,P^{2},P^{3}\}$.
\end{itemize}
\noindent {\bf Theorem 5.6.} Let $H$ and $N$ be normal subgroups
of a
 group $G$ and $X=(R;S)$
be a pseudo-Cayley graph. Then
\begin{itemize}
    \item[(1)]  $\overline{X}'_{H \cap N} = \overline{X}'_{H } \cap
    \overline{X}'_{N}$,
    \item[(2)]  $\underline{X}'_{H \cap N} = \underline{X}'_{H } \cap
    \underline{X}'_{N}$.
\end{itemize}
\noindent {\it Proof.}
\begin{itemize}
    \item[(1)] We have $$
    \begin{array}{ll}
    \overline{X}'_{H \cap N} & = ((H \cap
    N)^{\wedge}(R);S)\\
    & = (H^{\wedge}(R) \cap N^{\wedge}(R); S)\\
    &= (H^{\wedge}(R);S) \cap (N^{\wedge}(R);S)\\
    & =\overline{X}'_{H} \cap \overline{X}'_{N}.
    \end{array}
    $$
    \item[(2)]  By Theorem 3.2 (2), the proof is similar to (1).
\end{itemize}
\noindent{\bf Theorem 5.7.} Let $N$  be a normal subgroup, $H$ be
a subgroup of a group $G$ and $S \subseteq H$. Let $X=(H;S)$ be a
Cayley graph. Then

\begin{itemize}
    \item[(1)]  $N \subseteq H \Leftrightarrow \underline{X'} = X =
 \overline{X'}$,
    \item[(2)]  $N \not\subseteq H \Leftrightarrow \underline{X'} =
    \emptyset$.
\end{itemize}
{\it Proof.} \begin{itemize}
    \item[(1)]  In order to prove this term,  we show that $$N_{-}(H) =
 H =
    N^{\wedge}(H) \Longleftrightarrow N \subseteq H .$$

    Suppose that $N \subseteq H$. If $h \in H$, then $ hN \subseteq HN
 \subseteq HH (N
    \subseteq H) \subseteq H (H \leq G)$. So $h \in N_{-}(H)$ which
 implies that
    $H \subseteq N_{-}(H)$. Now, by Theorem 3.1 (1)  and the above
 result, we obtain  $H = N_{-}(H)$.
    If $g \in N^{\wedge}(H)$, then $gN \cap H \neq
    \emptyset$. So there exists $n \in N$ such that $gn \in H$ which
 implies that
    $gnn^{-1} \in H$ $(N \subseteq H, n$ and $n^{-1}
    \in H)$. Hence $g \in H$ and so $N^{\wedge}(H) \subseteq H$.
 Therefore, by Theorem 3.1 (1)  and the above results
    we have $N^{\wedge}(H) = H$.\\
    \indent Conversely $N_{-}(H) = H = N^{\wedge}(H)$ implies that
    $ 1N = N \subseteq H$.
 \item[(2)]  If $g \in N_{-}(H)$, then  $gN \subseteq H$. So for all $n
 \in N$, $gn \in H$. Since  $N_{-}(H) \subseteq H$, then $g \in H$.
 Thus
        $$g^{-1} \in H\Rightarrow g^{-1}gn \in H\Rightarrow n \in
 H\Rightarrow N \subseteq H $$

    and this shows repugnance.
\end{itemize}

Notice that according to Theorem 3.3, $N_{-}(H)$ and
    $N^{\wedge}(H)$ are subgroups of $G$. Then $\underline{X}$ and
    $\overline{X}$ are Cayley graphs and $X$ is definable. (A
    subset $X$ of $U$ is called {\it definable} if $\underline{Apr}(X)
 =
    \overline{Apr}(X)$).\\
\indent By Theorem 3.3,  $N^{\wedge}(H)$ is a subgroup of $G$.
Then
    $\overline{X}$ is a Cayley graph.\\ \\
{\bf Theorem 5.8.} Let $N$  be a normal subgroup, $R$ and $S$ be
subsets of a group $G$ and $SR \subseteq R$. Let $X=(R;S)$ be a
pseudo-Cayley graph. If $N_{-}(R)$ is not empty and there exists
$r \in R$ such that
 $ <s>r = R$, then $X$ is definable.\\ \\
{\it Proof.}  Since $N_{-}(R)$ is not empty, then  there exists
$r' \in N_{-}(R)$ and there exist $s_{1}, s_{2}, \ldots s_{m}$
such that $r' =
 s_{1}s_{2} \ldots
s_{m}r$. We have
$$rN = s_{1}s_{2} \ldots s_{m}r'N \subseteq s_{1}s_{2} \ldots
s_{m}R \subseteq R\ \  (SR \subseteq R ,\ r' \in N_{-}(R))$$ so $r
\in N_{-}(R)$. Thus for all $r'' \in R$ there exist $s'_{1}s'_{2}
 \ldots s'_{m'}$ such that $r'' =
s'_{1}s'_{2} \ldots s'_{m'}~r$. We have
$$r''N =  s'_{1}s'_{2} \ldots s'_{m'}~rN \subseteq s'_{1}s'_{2}
\ldots s'_{m'}~R \subseteq R \ \  (SR \subseteq R ,\ r \in
N_{-}(R))$$ which implies that $r'' \in N_{-}(R)$. Therefore we
obtain
$ N_{-}(R) = R = N^{\wedge}(R)$. Thus $X$ is definable.\\ \\
\indent Notice that the converse of Theorem 5.8 may  not be true.
For example, let  $G$ be the dihedral group of order $8$. Let $R =
\{P, P^{2}, P^{3} ,P \varepsilon,  P^{2} \varepsilon, P^{3}
\varepsilon
 \} ,$ $ S = \{\varepsilon \}$ and $N = \{1\}$. Then $X =
(R;S)$ is definable.
\\ \\
$(\underline{X}',\overline{X}')$ is called a  {\it rough vertex
pseudo-Cayley graph $X = (R;S)$}. A pseudo-Cayley graph $X =
(R;S)$ is called an {\it $N^{\wedge}$-vertex rough generating} if
$R$ is $N^{\wedge}$-rough subgroup of $G$ and  $S$ is a generating
set for $N^{\wedge} (R)$. Similarly, a pseudo-Cayley graph $X =
(R;S)$ is called an {\it $N_{-}$-vertex rough generating} if $R$
is $N_{-}$-rough subgroup of $G$ and $S$ is a generating set for
$N_{-}(R)$.\\ \\
A pseudo-Cayley graph $X = (R;S)$  is called an {\it
$N^{\wedge}$-vertex rough optimal connected} if  $R$ is
$N^{\wedge}$-rough subgroup of $G$ and  $S$ is a minimal Cayley
set for $N^{\wedge} (R)$. Similarly, a pseudo-Cayley graph $X =
(R;S)$ is called an {\it $N_{-}$-vertex rough optimal} if $R$ is
$N_{-}$-rough subgroup of $G$ and $S$ is a minimal Cayley set
for $N_{-}(R)$.\\ \\
{\bf Theorem 5.9.} Let $X=(R;S)$ be a pseudo-Cayley graph. If $X$
is an $N^{\wedge}$-vertex rough generating, then $\overline{X}'$
is connected. Similarly, if $X$ is an $N_{-}$-vertex rough
generating
then $\underline{X}'$ is connected.\\ \\
{\it Proof.} It is straightforward.\\ \\
{\bf Theorem 5.10.} Let $X=(R;S)$ be a pseudo-Cayley graph. If $X$
is
 an
$N^{\wedge}$-vertex rough optimal connected, then $\overline{X}'$
is optimal connected. Similarly, if $X$ is an
 $N_{-}$-vertex rough optimal connected, then $\underline{X}'$ is
optimal connected.\\ \\
{\it Proof.} It is straightforward.
\section{Rough pseudo-Cayley graphs}
\begin{figure}
\begin{center}
\includegraphics[1mm,0mm][85mm,30mm]{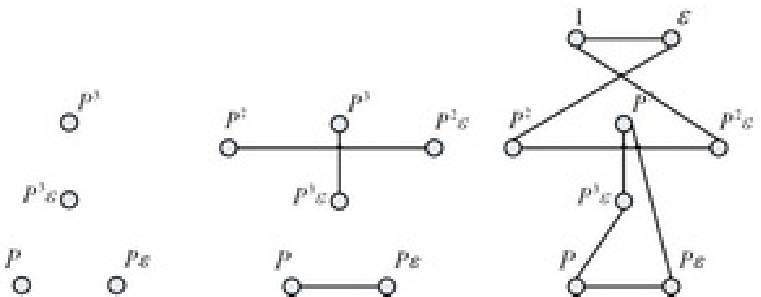}
\caption{The above  graphs are, respectively $\underline{X}''$,
$X$ and $\overline{X}''$.} \label{fig3}
\end{center}
\end{figure}

In this section, concept of lower and upper approximations
pseudo-Cayley graphs of a pseudo-Cayley graph with respect to a
normal subgroup is introduced then some
 properties of lower and upper approximations are brought.\\ \\
{\bf Definition 6.1.} Let $G$ be a finite group with identity $1$,
$N$
 be a normal
subgroup, $R$ be a subset of $G$, $X=(R;S)$ be a pseudo-Cayley
graph. Then the following graphs (we will prove these graphs are
 pseudo-Cayley
graphs):
$$\overline{X}'' = (N^{\wedge} (R);N^{\wedge} (S)^{\ast}) \ \ {\rm
  and} \ \ \underline{X}'' = (N_{-}(R);N_{-}(S))$$
are called, respectively,  {\it lower} and {\it upper
approximations pseudo-Cayley graphs} of a pseudo-Cayley graph
$X(R;S)$ with respect to the normal subgroup $N$.\\ \\
{\bf Theorem 6.2.} The two graphs $\underline{X}''$  and
$\overline{X}''$ are
pseudo-Cayley graphs.\\ \\
{\it Proof.} If $a \in N^{\wedge} (S)^{\ast}$ and $ b \in
N^{\wedge} (R)$, then $aN \cap S \neq \emptyset$ and $bN \cap R
\neq \emptyset$.  So there exist $n_{1},~n_{2} \in N$ such that
$an_{1} \in S$ and $bn_{2} \in R$. Hence $an_{1}bn_{2} \in R$
which implies that $an_{1}bn_{2} \in aNbN = abN$. We obtain $abN
\cap R \neq \emptyset $ which implies that $ab \in N^{\wedge}
(R)$. Therefore $N^{\wedge} (S)^{\ast}N^{\wedge} (R) \subseteq
N^{\wedge} (R)$.\\
\indent If $a \in N_{-}(S)$ and $ b \in N_{-}(R)$, then $aN
\subseteq S $ and $bN \subseteq R$. We have $$ab N = aNbN
\subseteq
 SR \subseteq R$$
so $ab \in N_{-}(R)$. Thus $N_{-}(S)N_{-}(R) \subseteq N_{-}(R)$.
Since
 $S \subseteq R$, then $N^{\wedge} (S)^{\ast} \subseteq
N^{\wedge} (S) \subseteq N^{\wedge} (R)$ and $N_{-}(S) \subseteq
N_{-}(R)$. Therefore $\underline{X}''$  and $\overline{X}''$ are
 pseudo-Cayley
graphs.\\ \\
{\bf Example 6.3.} Let $G$ be a dihedral group of order $8$. Let
$N = \{1, P^{2} \} $ be a normal subgroup and $R = \{P, P^{2},
P^{3} ,P \varepsilon , P^{2} \varepsilon , P^{3} \varepsilon  \}$
be a subset of $G$. Let pseudo-Cayley graph $X = (R;S)$ such that
$S$ equals to $\{\varepsilon \}$. We have
$$\underline{X}'' = (\{P, P^{3} ,P \varepsilon , P^{3} \varepsilon
\};\emptyset)$$ and $$ \overline{X}'' = (\{P, P^{2}, P^{3}, P
\varepsilon , P^{2} \varepsilon , P^{3} \varepsilon , 1,
\varepsilon \};\{\varepsilon, P^{2} \varepsilon\}).$$ (see figure
\ref{fig3}).
\\ \\
{\bf Theorem 6.4.} Let $N$ and $H$ be normal subgroups  and $R$,
$R_{1}$ and $R_{2}$ be subsets of a group $G$. Let $X=(R;S)$,
$X_{1}=(R_{1};S_{1})$ and $X_{2}=(R_{2};S_{2})$ be pseudo-Cayley
graphs. Then we have

\begin{itemize}
    \item [(1)] $\underline{X}'' \subseteq X \subseteq \overline{X}''$,

 \item[(2)]  $\underline{X_{1} \cap X_{2}}'' = \underline{X_{1}}'' \cap
\underline{X_{2}}'' $,

 \item[(3)]  $X_{1} \subseteq X_{2} \Longrightarrow \underline{X_{1}}''
 \subseteq \underline{X_{2}}''$,

 \item[(4)]  $X_{1} \subseteq X_{2} \Longrightarrow \overline{X_{1}}''
 \subseteq \overline{X_{2}}''$,

\item[(5)]  $\overline{X_{1} \cap X_{2}}'' \subseteq
\overline{X_{1}}'' \cap \overline{X_{2}}'' $,

\item[(6)]  $N \subseteq H \Rightarrow \overline{X_{N} }''
\subseteq \overline{X_{H}}''$,

\item[(7)]  $N \subseteq H \Rightarrow \underline{X_{H} }''
\subseteq \underline{X_{N}}''$.

\end{itemize}
{\it Proof.}
\begin{itemize}
    \item[(1)]  We have
    $$
    \begin{array}{ll}
    \underline{X}'' & = (N_{-}(R);N_{-}(S)) \subseteq
(R;N_{-}(S))\\
& \subseteq (R;S)\\
&\subseteq (N^{\wedge} (R);S)\\
& \subseteq (N^{\wedge} (R);N^{\wedge} (S)^{\ast})\\
&  =
    \overline{X}''.
\end{array}
$$
\item[(2)] We have
$$
\begin{array}{ll}
\underline{X_{1} \cap X_{2}}'' & = \underline{(R_{1}
\cap R_{2};S_{1} \cap S_{2})}''\\
&= (N^{\wedge} (R_{1} \cap R_{2});N^{\wedge} (S_{1} \cap
S_{2})^{\ast})\\
&= (N^{\wedge} (R_{1}) \cap N^{\wedge} (R_{2});N^{\wedge}
(S_{1})^{\ast} \cap N^{\wedge} (S_{2})^{\ast})\\
& = (N^{\wedge} (R_{1});N^{\wedge} (S_{1})^{\ast}) \cap
(N^{\wedge} (R_{2});N^{\wedge} (S_{2})^{\ast})\\
&= \underline{X_{1}}'' \cap \underline{X_{2}}'' .
\end{array}
$$

\item[(3)] Since $X_{1} \subseteq X_{2}$, then $R_{1 } \subseteq
R_{2}
 , S_{1 } \subseteq S_{2} $. Hence
$N_{-}(R_{1 }) \subseteq N_{-}(R_{2 }) , N_{-}(S_{1 }) \subseteq
N_{-}(S_{2 })$  which implies that $\underline{X_{1}}''
 \subseteq \underline{X_{2}}''$.

\item[(4)] By Theorem 3.1 (5),  the proof is similar to (3).

\item[(5)] We have
$$
\begin{array}{ll}
\overline{X_{1} \cap X_{2}}'' & = \overline{(R_{1} \cap
R_{2};S_{1} \cap S_{2})}''\\
&= (N^{\wedge}(R_{1} \cap R_{2});N^{\wedge}(S_{1} \cap
S_{2})^{\ast})\\
& \subseteq (N^{\wedge}(R_{1}) \cap
N^{\wedge}(R_{2});N^{\wedge}(S_{1})^{\ast} \cap
N^{\wedge}(S_{2})^{\ast})
\\
& = (N^{\wedge}(R_{1});N^{\wedge}(S_{1})^{\ast}) \cap
(N^{\wedge}(R_{2});N^{\wedge}(S_{2})^{\ast})\\
& = \overline{X_{1}}'' \cap \overline{X_{2}}''.
\end{array}
$$
\item[(6)]  Since $N \subseteq H$, then $N^{\wedge}(R) \subseteq
 H^{\wedge}(R)$ and $N^{\wedge}(S)
\subseteq H^{\wedge}(S)$. So
$$\overline{X_{N} }'' = (N^{\wedge}(R);N^{\wedge}(S))
\subseteq (H^{\wedge}(R);H^{\wedge}(S)) = \overline{X_{H}}''.$$

\item[(7)]  By Theorem 3.1 (9), the proof is similar to (8).
\end{itemize}
{\bf Remark 6.5.} Let $N$  be a normal subgroup and $H$ be a
subgroup of group $G$ such that $N \subseteq H$ and $S \subseteq
H$. Let $X=(H;S)$ be a Cayley graph. Then
 $(\underline{X}'',\overline{X}'')$  equals by an rough edge Cayley graph of $X$.\\ \\
{\it Proof.} This is simply provable by Theorem 5.7.\\ \\
\indent Notice that according to previous remark, and Examples
4.5, the
 converse part of the above items(4--7) are not
always true.\\ \\
{\bf Theorem 6.6.} Let $H$ and $N$ be normal subgroups of a group
$G$. Let $X=(R;S)$ be pseudo-Cayley graph. Then
\begin{itemize}
    \item[(1)]  $\overline{X}''_{H \cap N} = \overline{X}''_{H } \cap
    \overline{X}''_{N}$,
    \item[(2)]  $\underline{X}''_{H \cap N} = \underline{X}''_{H } \cap
    \underline{X}''_{N}$.
\end{itemize}
{\it Proof.}
\begin{itemize}
    \item[(1)] We have
    $$
    \begin{array}{ll}
    \overline{X}''_{H \cap N} & = ((H \cap
    N)^{\wedge}(R);(H \cap
    N)^{\wedge}(S)^{\ast})\\
    & = (H^{\wedge}(R) \cap N^{\wedge}(R); H^{\wedge}(S)^{\ast} \cap
 N^{\wedge}(S)^{\ast})\\
&= (H^{\wedge}(R);H^{\wedge}(S)^{\ast}) \bigcap
 (N^{\wedge}(R);N^{\wedge}(S)^{\ast})\\
& =\overline{X}''_{H } \bigcap \overline{X}''_{N}.
\end{array}
$$
    \item[(2)]  By Theorem 3.2 (2), the proof is similar to (1).
\end{itemize}
\noindent $(\underline{X}'',\overline{X}'')$ is called a  {\it
rough pseudo-Cayley graph $X = (R;S)$}. A pseudo-Cayley graph $X =
(R;S)$ is called an {\it $N^{\wedge}$-rough generating} if $R$ is
$N^{\wedge}$-rough subgroup of $G$ and $N^{\wedge} (S)$ is a
generating set for $N^{\wedge} (R)$. Similarly, a pseudo-Cayley
graph $X = (R;S)$ is called an {\it $N_{-}$-rough generating} if
$R$ is $N_{-}$-rough subgroup of $G$ and $N_{-}(S)$ is a
generating set for $N_{-}(R)$.
\\ \\
A pseudo-Cayley graph $X = (R;S)$  is called an {\it
$N^{\wedge}$-rough optimal connected} if  $R$ is
$N^{\wedge}$-rough subgroup of $G$ and $N^{\wedge} (S)$ is a
minimal Cayley set for $N^{\wedge} (R)$. Similarly, a
pseudo-Cayley graph $X = (R;S)$ is called an {\it $N_{-}$-rough
optimal} if $R$ is $N_{-}$-rough subgroup of $G$ and
the $N_{-}(S)$ is a minimal Cayley set for $N_{-}(R)$.\\ \\
{\bf Theorem 6.7.} Let $X=(R;S)$ be a pseudo-Cayley graph. If $X$
is an $N^{\wedge}$-rough generating, then $\overline{X}'' $ is
connected. Similarly, If $X$ is an $N_{-}$-rough generating, then
$\underline{X}''$ is connected.\\ \\
{\it Proof.}   According to Theorem 2.2, the proof is clear.\\ \\
{\bf Theorem 6.8.} Let $X=(R;S)$ be a pseudo-Cayley graph. If $X$
is an $N^{\wedge}$-rough optimal connected, then $\overline{X}''$
is optimal connected. Similarly, if $X$ is an
 $N_{-}$-rough optimal connected, then $\underline{X}'' $ is
optimal connected.\\
\\
{\it Proof.} According to Theorem 2.3, the proof is clear.

\section{Conclusion}
This paper  addressed a connection between two research fields,
rough set and Cayley graphs, which have a wide variety of
applications. Three approximations called rough edge Cayley
graphs, rough vertex pseudo-Cayley graphs and rough pseudo-Cayley
graphs on Cayley graphs and pseudo-Cayley graphs have been
defined. Some theorems and properties such as connectivity have
been discussed. In order to  find some examples and
contradictions, a  software is developed in C++ programming
language. These approximations can be applied to many challenging
problems in distributed systems such
as reliability and fault tolerance. The applications of these results are the purpose of a future study.\\
\\
{\bf Acknowledgement.} The authors are highly grateful to
referees for their valuable comments and suggestions for improving
the paper.

\end{document}